\documentclass[conference]{IEEEtran}
\IEEEoverridecommandlockouts
\usepackage{cite}
\usepackage{amsmath,amssymb,amsfonts}
\usepackage{algorithmic}
\usepackage{graphicx}
\usepackage{textcomp}
\usepackage{xcolor}
\usepackage[normalem]{ulem}

\def\BibTeX{{\rm B\kern-.05em{\sc i\kern-.025em b}\kern-.08em
    T\kern-.1667em\lower.7ex\hbox{E}\kern-.125emX}}
\begin{document}
\renewcommand{\slash}{/\penalty\exhyphenpenalty\hspace{0pt}}
\newenvironment{ldescription}[1]
  {\begin{list}{}%
   {\renewcommand\makelabel[1]{##1\hfill}%
   \settowidth\labelwidth{\makelabel{#1}}%
   \setlength\leftmargin{\labelwidth}
   \addtolength\leftmargin{\labelsep}}}
  {\end{list}}

\title{Pp OPF - Pandapower Implementation of Three-phase Optimal Power Flow Model
\thanks{This work has been supported in part by the European Structural and Investment Funds under KK.01.2.1.02.0042 DINGO (Distribution Grid Optimization) and in part by the FLEXIGRID project from the European Union's Horizon 2020 research and innovation programme under grant agreement No 864579. This paper reflects the FLEXIGRID consortium view and the European Commission is not responsible for any use that may be made of the information it contains.}
}

\author{\IEEEauthorblockN{Tomislav Antić}
\IEEEauthorblockA{\textit{Faculty of Electrical Engineering}\\\textit{and Computing} \\
\textit{University of Zagreb}\\
Zagreb, Croatia \\
tomislav.antic@fer.hr}
\and
\IEEEauthorblockN{Andrew Keane}
\IEEEauthorblockA{\textit{UCD Energy Institute} \\
\textit{University College Dublin}\\
Dublin, Ireland \\
andrew.keane@ucd.ie}
\and
\IEEEauthorblockN{Tomislav Capuder}
\IEEEauthorblockA{\textit{Faculty of Electrical Engineering}\\\textit{and Computing} \\
\textit{University of Zagreb}\\
Zagreb, Croatia \\
tomislav.capuder\@fer.hr}
}

\maketitle

\begin{abstract}
Challenges in the planning and operation of distribution networks caused by the integration of distributed energy resources (DERs) create the need for the development of tools that can be easily used by system operators, industry, and the research society but are also easily upgraded with new functionalities. The full implementation of one such open source tool, named \textit{pp OPF} (pandapower Optimal Power Flow), is presented in this paper. \textit{Pp OPF} is the tool used for three-phase optimal power flow (OPF) calculations and it is based on already existing functionalities of pandapower, a Python library for power system calculations. The developed tool enables the possibility to use both power-voltage and current-voltage formulation in different OPF problems, such as determining the photovoltaic (PV) hosting capacity in three-phase distribution networks, the problem on which the functionality of the developed tool is tested. Additionally, the accuracy of \textit{pp OPF} is verified by comparing the results with ones obtained by pandapower power flow calculation for the same set of input values. The open-source implementation allows further upgrades, the addition of new functionalities, and the creation of new case studies relevant to the planning and operation of distribution networks.
\end{abstract}

\begin{IEEEkeywords}
distribution networks, distributed energy resources, pandapower, three-phase optimal power flow
\end{IEEEkeywords}

\section*{Nomenclature} 

\subsection{Sets and Indices}
\begin{ldescription}{$xxxxxx$}
    \item[$N$] Set of all nodes
    \item[$P$] Set of all phases
    \item[$T$] Set of all time intervals
    \item[$i, j$] Begin and end node of the element, $i, j \in N$
    \item[$n$] Observed node, $n \in N$
    \item[$slack$] Slack node
    \item[$p$] Phase of the node, $p \in P$
    \item[$t$] $t \in T$
\end{ldescription}
\subsection{Parameters}
\begin{ldescription}{$xxcxxxxxxxxxxx$}
    \item[$\begin{bmatrix} Z_{012} \end{bmatrix}$] Sequence impedance matrix
    \item[$\begin{bmatrix} A \end{bmatrix}$] Transformation matrix
    \item[$P_{i,p,t \ l}$] Active power of load $l$ at phase $p$ and in time period $t$
    \item[$Q_{i,p,t \ l}$] Reactive power of load $l$ at phase $p$ and in time period $t$
    \item[$U^{min/max}$] Minimum and maximum voltage magnitude
    \item[$VUF^{max}$] Maximum value of voltage unbalance factor
    \item[$S_{max,ij}$] Maximum apparent power of element $ij$
\end{ldescription}
\subsection{Variables}
\begin{ldescription}{$xxxxxxxxx$}
    \item[$U_{n,p,t}^{re/im}$] Real/imaginary part of voltage at node $n$, phase $p$, and in time period $t$
    \item[$P_{ij,p,t}$] Active power between nodes $i$ and $j$, of phase $p$, and in time period $t$
    \item[$Q_{ij,p,t}$] Reactive power between nodes $i$ and $j$, of phase $p$, and in time period $t$
    \item[$P_{n,p,t \ l/g}$] Active power of load $l$/generator $g$  connected to the node $n$, phase $p$, in time period $t$
    \item[$Q_{n,p,t \ l/g}$] Reactive power of load $l$/generator $g$  connected to the node $n$, phase $p$, in time period $t$
\end{ldescription}

\section{Introduction}
The share of low-carbon (LC) technologies has been continuously increasing. Despite the positive environmental and financial effects, the integration of LC units is often uncoordinated which can lead to different technical problems in the operation of distribution networks, including overvoltage and undervoltage \cite{9808683}, increased voltage unbalance \cite{en14010117}, and overloading of lines and transformers \cite{9935452}.

Despite the high number of commercial software used for distribution networks analyses, recent digitization and more often exploitation of computing techniques in power systems planning and operation has led to the development of open source tools in high-level programming languages, such as MATPOWER \cite{5491276} or pandapower \cite{8344496}. These and similar tools are often used only for calculations of technical conditions in distribution networks and in most cases do not have implemented optimization techniques that can help in decision-making. As a way of battling these problems, the distribution system operators now more than ever reach for optimization tools as opposed to the conventional simulation-based approach, making optimal power flow algorithms essential in everyday operation.

Even though OPF models are well-known and widely used in different power system problems, e.g., unit commitment \cite{9305983}, their use in distribution and especially low voltage (LV) networks require modification and transition towards three-phase four-wire models. Three-phase OPF models have numerous applications, including optimal control of devices \cite{8457301}, congestion mitigation \cite{VANIN2022100936}, calculation of operating envelopes \cite{9816082}, etc. Three-phase optimal power flow problems can be formulated as current-voltage \cite{9133699} or power-voltage \cite{geth_power_voltage} and can be modeled as exact, non-convex problems \cite{9133699} or they can be linearized \cite{9816082} or relaxed \cite{HUANG2022117771}.

Recent developments of open source software include the implementation of a three-phase OPF model as a Julia package \cite{FOBES2020106664} and integration of a three-phase OPF model with OpenDSS software \cite{9282125}. Despite the recent developments in pandapower, including the creation of new elements and introduction of new functionalities, e.g., three-phase power flow and state estimation \cite{9769084} or the availability of harmonic power flow simulations \cite{ANTIC2022107935}, using pandapower Python library in three-phase OPF problems is still not possible. 

To overcome the detected pandapower shortfall, we propose an implementation of non-linear, current-voltage, and power-voltage three-phase OPF formulations using existing pandapower and newly added functionalities necessary for modeling the wanted problem. The developed tool named \textit{pp OPF} (pandapower Optimal Power Flow) is open source and can be used by accessing \cite{ppopf}. The functionality of the tool is tested on the problem of calculating photovoltaic (PV) hosting capacity in a real-world LV network. The accuracy of the tool is verified against the results calculated by pandapower power flow simulations for the same input data. It is important to emphasize that the presented mathematical model and defined case study on which the model was tested are not novel contributions. However, to the best of the authors' knowledge, this is the first time that pandapower was used for the implementation of a three-phase OPF mathematical model. Additionally, problems that are solved by exact non-linear models are often non-scalable and require a significant computational burden in order to be solved. In this paper, an introduction to such problems and one possible method for overcoming them is presented next to the verification and the analysis of the case study results. \textit{Pp OPF} is the only three-phase OPF tool that has been built upon the widely used Python library pandapower, and together with the development of Open-DSOPF presented in \cite{9282125}, it is the only publicly available Python-based tool that enables OPF simulations in three-phase distribution networks. As shown later in the paper, \textit{pp OPF} is accurate, efficient, and provides a feasible solution, the same as OpenDS-OPF. Therefore, the selection of the tool for solving the given problems completely depends on the preferences of a user, but the presented development gives another option to researchers and the industry.

The rest of the paper is organized as follows: Section \ref{sec:pp_connection} presents the connection to pandapower by using some of the already developed functionalities. In Section \ref{sec:opf_model}, the mathematical model of the power-voltage formulation is given. The results of a defined case study and the comparison of the results with ones obtained by pandapower are presented in Section \ref{sec:cs_ver}, while the conclusions and future work are given in Section \ref{sec:conclusions}.

\section{Connection to pandapower}\label{sec:pp_connection}
Pandapower is developed as an open source tool for power system modeling, analysis, and optimization with a high degree of automation. Since the new tool presented in this paper is used for analyses and optimization problems that have not yet been implemented in pandapower, the only existing functionality that is used is the modeling of distribution networks by using pandapower for defining different elements necessary for creating the representation of the observed distribution network, while all other necessary functions are built on top of that. Used elements are an external grid, transformers, lines, and loads. An external grid is used to define the per-unit voltage at the slack node. The assumption in the presented three-phase OPF model is that the voltage of the slack node is symmetrically distributed over phases, as shown in eq. (\ref{eq:slack_bus_voltage}).

\begin{equation}
\begin{gathered}
    U_{slack,a} = U_{slack} \\
    U_{slack,b} = U_{slack}\angle-120^{\circ} \\
    U_{slack,c} = U_{slack}\angle120^{\circ}
\end{gathered}
\label{eq:slack_bus_voltage}
\end{equation}

Transformers are defined with parameters needed for the calculation of their impedance, such as voltage ratio, rated power, and short-circuit voltage. Lines are created with parameters including but not limited to resistance and reactance for zero and positive sequence systems, length, and maximum current. After creating these elements, already existing pandapower functionalities are used for calculating their impedance and placing them in the impedance matrix $\begin{bmatrix}Z_{012}\end{bmatrix}$ containing the calculated impedance of zero, positive, and negative sequence systems. By using transformation matrix $\begin{bmatrix}A\end{bmatrix}$ and its inverse $\begin{bmatrix}A\end{bmatrix}^{-1}$, three-phase impedance matrix $\begin{bmatrix}Z_{abc}\end{bmatrix}$ is calculated.

All loads in a distribution network are created using pandapower functionalities, together with other elements relevant to the representation of the network's mathematical mode. Since \textit{pp OPF} is used in the multi-temporal optimization of three-phase power systems, loads need to be defined with the active and reactive power of each phase in each time period. Since the use of pandapower in multi-temporal simulations of a three-phase network is not straightforward, a load curve that defines the change of the loads' value during the observed period is additional input in the tool. Based on the objective function and the goal of the optimization problem, additional elements representing DERs such as PVs or electric vehicles (EVs) can be added to the formulation. These elements can be defined with a fixed value or can be defined as variables whose power needs to be determined during the solving of the defined problem.

In order to integrate the presented tool with already existing pandapower functionalities, a Python script that creates all of the network's elements needs to be created first. After the creation of the pandapower model of a three-phase network, a script modifies the network data and creates three-phase impedance and admittance matrices. Finally, all mathematical equations and the constraints valid for both current-voltage and power voltage formulations of the three-phase OPF model are integrated into two final scripts. A user of the tool is only responsible for the creation of the mathematical model of a network, in a same way as in the case of using pandapower, i.e., there is no need for additional effort since the functionalities of the OPF model are developed as a part of the source code.

\section{Three-phase OPF mathematical model}\label{sec:opf_model}
After defining a three-phase distribution network with all of its elements using pandapower, its parameters are used in the developed OPF formulation. The \textit{pp OPF} tool integrates both power-voltage and current-voltage formulations. OPF models are formulated using Python programming language and Pyomo optimization framework. However, in this paper, only the mathematical model of a power-voltage formulation is presented, while both formulations are tested and verified on case studies presented in Section \ref{sec:cs_ver}.

Active and reactive power flow between nodes $i$ and $j$ is described with eqs. (\ref{eq:line_act_power_calc_det}) and (\ref{eq:line_react_power_calc_det}). Additionally, the power flow through each element cannot be larger than the maximum power which is defined as an input parameter of a transformer or a line in the observed network. This constraint is given with eq. (\ref{eq:max_power}).

\begin{equation}
\begin{gathered}
    P_{ij,p,t} = \sum_{q \in\{a,b,c\}}(U_{i,p,t}^{re}\cdot U_{i,q,t}^{re} + U_{i,p,t}^{im}\cdot U_{i,q,t}^{im})\cdot G_{ij,pq} + \\
    \sum_{q \in\{a,b,c\}}(U_{i,p,t}^{im}\cdot U_{i,q,t}^{re} - U_{i,p,t}^{re}\cdot U_{i,q,t}^{im})\cdot B_{ij,pq} -\\
    \sum_{q \in\{a,b,c\}}(U_{i,p,t}^{re}\cdot U_{j,q,t}^{re} + U_{i,p,t}^{im}\cdot U_{j,q,t}^{im})\cdot G_{ij,pq} -\\
    \sum_{q \in\{a,b,c\}}(U_{i,p,t}^{im}\cdot U_{j,q,t}^{re} - U_{i,p,t}^{re}\cdot U_{j,q,t}^{im})\cdot B_{ij,pq}
\end{gathered}
\label{eq:line_act_power_calc_det}
\end{equation}

\begin{equation}
\begin{gathered}
    Q_{ij,p,t} = -\sum_{q \in\{a,b,c\}}(U_{i,p,t}^{re}\cdot U_{i,q,t}^{re} + U_{i,p,t}^{im}\cdot U_{i,,tq}^{im})\cdot B_{ij,pq} + \\
    \sum_{q \in\{a,b,c\}}(U_{i,p,t}^{im}\cdot U_{i,q,t}^{re} - U_{i,p,t}^{re}\cdot U_{i,q,t}^{im})\cdot G_{ij,pq} +\\
    \sum_{q \in\{a,b,c\}}(U_{i,p,t}^{re}\cdot U_{j,q,t}^{re} + U_{i,p,t}^{im}\cdot U_{j,q,t}^{im})\cdot B_{ij,pq} -\\
    \sum_{q \in\{a,b,c\}}(U_{i,p,t}^{im}\cdot U_{j,q,t}^{re} - U_{i,p,t}^{re}\cdot U_{j,q,t}^{im})\cdot G_{ij,pq}
\end{gathered}
\label{eq:line_react_power_calc_det}
\end{equation}

\begin{equation}
    P_{ij,p}^{re \ 2} + Q_{ij,p}^{im \ 2} \leq S_{max,ij}
\label{eq:max_power}
\end{equation}

Besides power flow, the only other constraint needed in the power-voltage formulation is Kirchoff's Current Law (KCL) presented with active and reactive power instead of current. For securing KCL, the following eqs. (\ref{eq:kirchoff_active}) and (\ref{eq:kirchoff_reactive}) are defined.

\begin{equation}
    P_{i,p,t \ g} + P_{h \rightarrow i,p,t} = P_{i,p,t \ l} + P_{i \rightarrow j,p,t}
    \label{eq:kirchoff_active}
\end{equation}

\begin{equation}
    Q_{i,p,t \ g} + Q_{h \rightarrow i,p,t} = Q_{i,p,t \ l} + Q_{i \rightarrow j,p,t}
    \label{eq:kirchoff_reactive}
\end{equation}

Equations used for constraining values of voltage magnitude and voltage unbalance factor (VUF) are presented with eqs. (\ref{eq:node_volt_limit}) and (\ref{eq:vuf}).

\begin{equation}
    (U^{min})^2 \leq (U_{i,p,t}^{re})^2 + (U_{i,p,t}^{im})^2 \leq (U^{max})^2
    \label{eq:node_volt_limit}
\end{equation}

\begin{equation}
   \frac{|U_{n,a} + a^2\cdot U_{b,n} + a\cdot U_{c,n}|^2}{|U_{n,a} + a\cdot U_{b,n} + a^2\cdot U_{c,n}|^2} \leq (VUF^{max})^2
    \label{eq:vuf}
\end{equation}

\section{PVs hosting capacity}\label{sec:cs_ver}
PVs or other DERs in general can be single-phase or three-phase connected to the network. Single-phase connection is easier and cheaper at the LV level but is limited with relatively low export or import power, e.g., export power in Croatian LV networks may not be higher than 3.68 kW \cite{CroatianGridCode}, and can lead to the increase in voltage unbalance. On the contrary, a three-phase connection allows the installation of PVs with higher production power and solves the unbalance problem but the connection is more expensive and complex.

In this paper, two case studies are defined for calculating the hosting capacity of PVs (generators $g$) in a real-world LV distribution network. In the first case study, PVs can only be single-phase connected to one phase at a time, and in the second, PVs are three-phase connected to the network. To obtain an exact, optimal solution, the presented formulation needs to be extended with binary variables that will define the exact connection phase in the case of a single-phase connection or if PVs are single-phase or three-phase connected at a certain node. Such an approach would uplift the formulation from the non-linear programming (NLP) problem to the mixed-integer (MI)NLP problem, which would lead to problems related to the scalability and the computational time needed for finding an optimal solution. Therefore, both case studies are defined as NLP problems and solved using the pyomo optimization framework and the ipopt solver. In the case of a single-phase connection of PVs, their connection phase is randomly pre-defined in each node. When a three-phase connection is observed, the power of PVs is distributed over the phases and PVs can not be single-phase connected. In cases of well-built and resilient networks, determining the hosting capacity will mostly lead to the connection of PVs with the larger connection power, which makes the assumption that PVs are always three-phase connected valid in solving the presented problem. In both observed case studies, the production power of PVs is constrained with the curve containing per-unit values scaled according to the largest measured value of production in the observed network's location, as shown in eq. (\ref{eq:pv_curve}).

\begin{equation}
    P_{n,p,t \ PV} \leq  P_{PV,max} \cdot P_{PV,t}
    \label{eq:pv_curve}
\end{equation}

where $P_{n,p,t \ PV}$ is the variable representing the calculated power of PV at the certain node $n$, phase $p$, and time period $t$, $P_{PV,max}$ is the theoretical maximum value of PVs production, and $P_{PV,t}$ is the per-unit value that ensures the value of production that does not exceed the limit defined by the irradiation and other physical quantities.

\begin{figure}[htbp]
    \centering
    \includegraphics[width=\columnwidth]{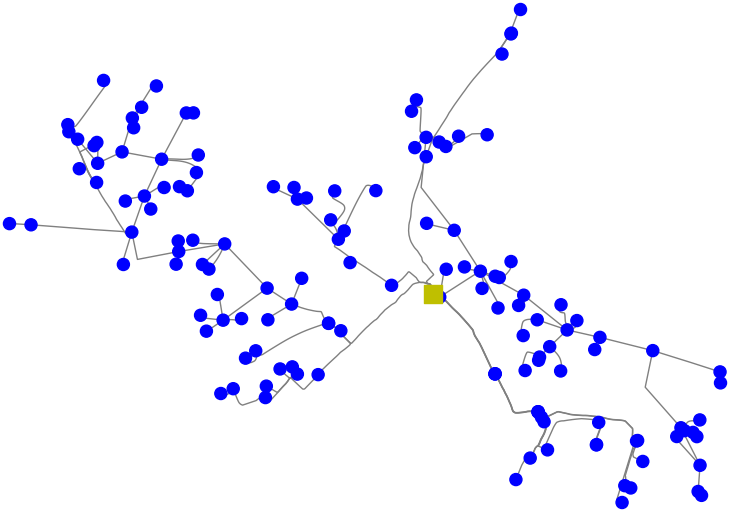}
    \caption{A residential LV network}
    \label{fig:net}
\end{figure}

For the presented case study, a real-world Croatian LV network shown in Fig. \ref{fig:net} was created using pandapower. Additionally, measurements collected from the smart meters were used to define the values of the loads' active and reactive power. 

\begin{figure}[htbp]
    \centering
    \includegraphics[width=\columnwidth]{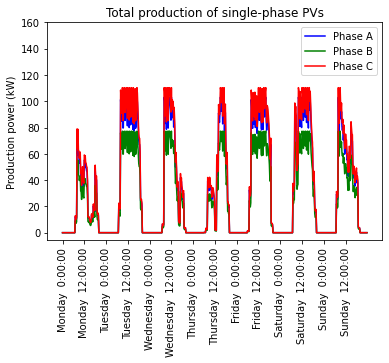}
    \caption{Production of single-phase PVs}
    \label{fig:pv_single_result}
\end{figure}

Fig. \ref{fig:pv_single_result} shows the summarized PV production curve for each phase over time in the first case study. The curve has the shape characteristic for the production of PVs, with the characteristic distortions caused by sudden weather changes that affect the production power. Together with the shape, values in the PV production curve are similar every day, with a maximum power of almost $120 kW$. The calculated value of the objective function is equal to the theoretical maximum, i.e., a different selection of PV connection phases would not lead to additional generation and the use of the NLP instead of MINLP formulation does not cause the loss of accuracy. However, a different network topology or a change of the objective function could change that and show the need for using the exact MINLP formulation. A more detailed investigation of that problem is out of the scope of the paper.

\begin{figure}[htbp]
    \centering
    \includegraphics[width=\columnwidth]{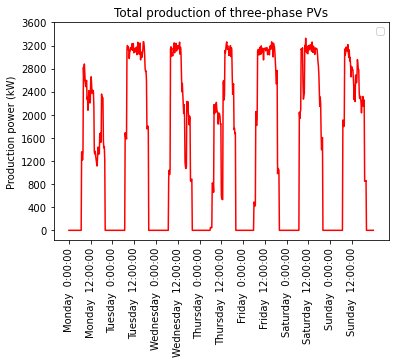}
    \caption{Production of three-phase PVs}
    \label{fig:pv_three_result}
\end{figure}

In the second case, PVs were three-phase connected to the network. Fig. \ref{fig:pv_three_result} does not show the power of each phase but only cumulative power at each node, since the production power at every phase is almost the same. The results show that the three-phase connection allows the significantly higher capacity of installed PVs in the same network, with the production power at several time periods of almost $3500 kW$. The results also confirm and make valid the assumption that none of the PVs are single-phase connected, which is relevant in avoiding MINLP formulations of a given problem.

The results clearly indicate that in order to maximize the share of LC units, a three-phase connection should be encouraged but a detailed study on finding the optimal location and size of PVs in the network is outside the scope of the paper. Moreover, it is already well-investigated in the existing literature. The main goal of this paper was to test the functionality of the implemented tool on one possible example. Additionally, the presented case study is used as an introduction to the potential methods that can be used in the improvement of the scalability of the problems solved using three-phase OPF formulations. 

Since both current-voltage and power-voltage formulations are implemented from scratch, it is necessary to verify the results of OPF calculations. Verification is made by comparing the results obtained by the \textit{pp OPF} with the results obtained by pandapower PF calculation for the same network and set of values. One of the outputs of the presented tool \textit{pp OPF} is the power of a single-phase or three-phase generator connected to a node. Together with the technical parameters of an LV network and active and reactive power of end-users' demand, the power of generators is used as input in the pandapower power flow calculation. Therefore, voltage magnitude values are calculated by two different tools for the same conditions, which enables the comparison of the results. Table \ref{tab:pp_ppopf} shows the comparison of voltage magnitude calculated with \textit{pp OPF} and pandapower. 

\begin{table}[htbp]
\caption{Comparison of results: pandapower and \textit{pp OPF}}
\label{tab:pp_ppopf}
\resizebox{\columnwidth}{!}{\begin{tabular}{c|c|c|c|c|c|}
\cline{2-6}
                                                                                                     & \begin{tabular}[c]{@{}c@{}}Maximum\\ error (p.u.)\end{tabular} & \begin{tabular}[c]{@{}c@{}}Minimum\\ error (p.u.)\end{tabular} & \begin{tabular}[c]{@{}c@{}}Average\\ error (p.u.)\end{tabular} & \begin{tabular}[c]{@{}c@{}}Median\\ error (p.u.)\end{tabular} & \begin{tabular}[c]{@{}c@{}}RMSE\\ (p.u.)\end{tabular} \\ \hline
\multicolumn{1}{|c|}{\textbf{\begin{tabular}[c]{@{}c@{}}Current-voltage\\ formulation\end{tabular}}} & 0.0624                                                         & 0.0000                                                         & 0.0024                                                         & 0.0000                                                        & 0.0086                                                \\ \hline
\multicolumn{1}{|c|}{\textbf{\begin{tabular}[c]{@{}c@{}}Power-voltage\\ formulation\end{tabular}}}   & 0.0472                                                         & 0.0000                                                         & 0.0074                                                         & 0.0055                                                        & 0.0097                                                \\ \hline
\end{tabular}}
\end{table}

As it can be seen from the comparison, neither of the calculated errors is of a significant value, with the maximum error lower than 0.07 p.u., which can be characterized as an outlier value and with even lower values of other calculated deviations. Calculated differences clearly indicate that the presented tool can be used in different problems in three-phase distribution networks.

\section{Conclusions and future work}\label{sec:conclusions}
In this paper, we present the development and implementation of \textit{pp OPF}, a three-phase OPF-based tool presented in this paper. \textit{pp OPF} is developed as an extension of pandapower, an already developed open source tool for power system simulations. Both power-voltage and current-voltage OPF formulations were added on top of existing pandapower functionalities.

Following the development and implementation, \textit{pp OPF} has been verified and tested on the real-world case study. In the defined case study, a real-world three-phase LV network was used in order to calculate its PV hosting capacity, both in cases of single-phase and three-phase connection of PVs. The results of the simulation were expected and emphasize the need for the encouragement of a three-phase connection or relaxing the limitations of a single-phase connection. A more detailed analysis of the results is outside the scope of this paper since the main objective was to present one of the potential functionalities of the developed tool.

The results of the comparison show that there are no significant differences between voltages calculated by the two different tools, with the maximum error detected as an outlier value and together with the values of other calculated errors, it clearly shows that there are no obstacles in using the \textit{pp OPF} tool in calculating PV hosting capacity and other optimization problems in three-phase distribution networks.

\bibliographystyle{IEEEtran}
\bibliography{IEEEabrv,paper_final}

\end{document}